\documentclass{article}
\usepackage[psamsfonts]{amssymb}
\newtheorem{theorem}{Theorem}[section]
\newtheorem{lemma}[theorem]{Lemma}
\newtheorem{e-proposition}[theorem]{Proposition}

\newtheorem{e-definition}[theorem]{Definition\rm}

\title{Asymptotic law of likelihood ratio for multilayer perceptron models}

\author{
Joseph Rynkiewicz\\
SAMM\\
Universit\'e de Paris 1\\
Paris, France\\
\texttt{joseph.rynkiewicz@univ-paris1.fr} \\
}

\begin{document}

\maketitle

\begin{abstract}
We consider regression models involving multilayer perceptrons (MLP) with one hidden layer and a Gaussian noise. The data are assumed to be generated by a true MLP model and the estimation of the parameters of the MLP is done by maximizing the likelihood of the model. When the number of hidden units of the true model is known,  the asymptotic distribution of the maximum likelihood estimator (MLE) and the likelihood ratio (LR) statistic is easy to compute and converge to a $\chi^2$ law. However, if the number of hidden unit is over-estimated the Fischer information matrix of the model is singular and the asymptotic behavior of the MLE is unknown. This paper deals with this case, and gives the exact asymptotic law of the LR statistics. Namely, if the parameters of the MLP lie in a suitable compact set, we show that the LR statistics is the supremum of the square of a Gaussian process indexed by a class of limit score functions.
\end{abstract}

\section{Introduction}
Feedforward neural networks are well known and are popular tools to deal with non-linear statistic models.
We can describe MLP as a parametric family of probability density functions. If the noise of the regression model is Gaussian then it is well known that the maximum likelihood estimator is equal to the least-square estimator. Therefore, Gaussian likelihood is the usual assumption when we consider feedforward neural networks from a statistical viewpoint. 
White [9] reviews statistical properties of  MLP estimation in detail. However he leaves an important question pending: the asymptotic behavior of the estimator when an MLP in use has redundant hidden units and the Fisher information matrix is singular. Amari, Park and Ozeki [1] give several examples of behavior of the LR in such cases. Fukumizu [4] shows that, for unbounded parameters, the  LR statistic can have an order lower bounded by $O(\log(n))$ with $n$ the number of observations instead of the classical convergence property to $\chi^2$ law. 

However, a fairly natural assumption is to consider that the parameters are bounded. Indeed, computer calculations  always assume that numbers are bounded. Moreover a safe practice is to bound the parameters in order to avoid numerical problems. In such context, different situations can occur. In some cases, such as mixture models, the LR is tight and the calculation of the asymptotic distribution is possible (see Liu and Shao [7]). In other cases it may occur that even if the parameters are bounded the likelihood ratio diverges this is for example the case in  hidden Markov models (see Gassiat and Keribin [5]). So the behavior of likelihood ratio in the case of MLPs with bounded parameters is still an open question. 

In this paper, we derive the distribution of the likelihood ratio if the parameters are in a suitable compact set (i.e. bounded and closed). To obtain this result we use recent techniques introduced by Dacunha-Castelle and Gassiat [2] and Liu and Shao [7]. These techniques consist in finding a parameterization separating the identifiable part and the unidentifiable part of the parameter vector, then we can obtain an asymptotic development of the likelihood of the model which allows us to show that a set of generalized score functions is a Donsker class and to find the asymptotic distribution of the LR statistic. The paper is organized as follows. In section 2 we state the model and the main assumptions. Section 3 presents our main theorem and explains its meaning with a brief summary and a statement of significance of this work. In section 4 we applied this theorem to the identification to the true architecture of the MLP function. In section 5 we show that MLP functions with sigmoidale transfert functions verify the assumption of this theorem. Finally, we prove the theorem in the appendix.

\section{The model}
We consider the model of regression for $i\in {\mathbb N}^*$:
\begin{equation}\label{model}
Y_i=F_{\theta^0}(X_i)+\varepsilon_i
\end{equation}
where $X_i\in\mathbb R^d$ are observed exogenous variables and $Y_i$ is the variable to explain. The data $(Y_i,X_i)$ are assumed to be generated by this true model. The noise $\left(\varepsilon_i\right)_{i\in {\mathbb N}^*}$ is a sequence of independent and identically distributed (i.i.d.) ${\cal N}(0,\sigma^2)$ variables. 
\subsection{The regression function}
Let  $x=(1,x_1,\cdots,x_d)^T\in{\mathbb R}^{d+1}$ be the vector of inputs and $w_i:=\left(w_{i0},w_{i1},\cdots,w_{id}\right)^T$, the MLP function with $k$ hidden units can be written : 

\[
F_\theta(x)=\beta+\sum_{i=1}^k a_i\phi\left(w_i^Tx\right),
\]
with $\theta=\left(\beta,a_1,\cdots,a_k,w_{10},\cdots,w_{1d},\cdots,w_{k0},\cdots,w_{kd}\right)\subset {\mathbb R}^{k\times (d+2)+1}$ the parameters of the model. The transfer function $\phi$ will be assumed bounded and three times derivable. We assume also that the first, second and third derivatives of the transfer function $\phi$:  $\phi^{'}$, $\phi^{''}$ and $\phi^{'''}$ are bounded. 
In order to simplify the presentation, we assume that the variance of the noise $\sigma^2$ is known. Note that it is assumed that the true model (\ref{model}) is included in the considered set of parameter $\Theta$. Let us define the true number of hidden units as the smallest integer $k^0$ so that $\theta^0=\left(\beta^0,a^0_1,\cdots,a^0_k,w^0_{10},\cdots,w^0_{1d},\cdots,w^0_{k0},\cdots,w^0_{kd}\right)$ exists with $F_{\theta^0}$ equal to the true regression function of model (\ref{model}). 
\subsection{Parameterization of the model}
Let us write $\Vert.\Vert$ for the Euclidean norm. Let us consider the variable  $Z_i=(X_i,Y_i)$ where $X_i$ and $Y_i$ follow the probability law induced by the model (\ref{model}). We assume that the law of $X_i$ will be $q(x)\lambda_d(x)$ with $\lambda_d$ the Lebesgue measure on ${\mathbb R}^d$ and $q(x)>0$ for all $x\in{\mathbb R}^d$. The likelihood of the observation $z:=(x,y)$ for a parameter vector $\theta=\left(\beta,a_1,\cdots,a_{k},b_1,\cdots,b_{k},w_{11},\cdots,w_{1d},\cdots,w_{{k}d}\right)$ will be written: 
\[
f_{\theta}(z)=\frac{1}{\sqrt{2\pi\sigma^2}}e^{-\frac{1}{2\sigma^2}\left(y-F_\theta(x)\right)^2}q(x).
\]
Let $\eta>0$ be a small constant and $M$ a huge constant, the set of possible parameters will be
\[
\begin{array}{l}
\Theta_k:=\left\{\theta=\left(\beta,a_1,\cdots,a_k,w_{10},\cdots,w_{1d},\cdots,w_{k0},\cdots,w_{kd}\right),\right.\\
\left.\forall 1\leq i\leq k, \Vert w_i\Vert\geq \eta,\Vert a_i\Vert\geq \eta\mbox{ and }\Vert\theta\Vert\leq M\right\}.
\end{array}
\]
\paragraph{Constraints on the parameter set.}

The constraint $\Vert w_i \Vert\geq \eta$ is introduced in order to avoid the hidden unit from being constant like the bias $\beta$, instead of being a function of $x$. The constraint $\Vert a_i\Vert\geq \eta$ forces the parameters of the hidden units to converge to one of the parameter vector $w^0_j, j\in\{1,\cdots,k^0\}$ when they maximize the likelihood. Finally, with the constraint $\Vert\theta\Vert\leq M$, the parameters are bounded and the set $\Theta_k$ compact. Note that these constraints are very easy to set in practice. 

The true density of the observation will be denoted $f(z):=f_{\theta^0}(z)$. The main goal of the parametric statistic is to give an estimation of the true parameter $\theta_0$ thanks to the observations $\left(z_1,\cdots,z_n\right)$. This can be done by maximizing the log-likelihood function :
\[
l_n(\theta):=\sum_{i=1}^{n}\log f_{\theta}(z_i).
\]
The parameter vectors ${\hat \theta}_n$ realizing the maximum will be called Maximum Likelihood Estimator (MLE). However,the MLE belongs to a non-null dimension submanifold if the number of hidden units is overestimated. In the next section we will study the behavior of 
\[
\sup_{\theta\in \Theta_k}\sum_{i=1}^{n}\log f_{\theta}(z_i)-\log f(z_i),\mbox{ where } k\geq k^0 
\]
which is the key to guess the true architecture of the MLP model.   
\section{Asymptotic distribution of the LR statistic}
We will use the abbreviation $P g=\int g dP$ for an integrable function $g$ and a measure $P$. We will define the $L^2(P)$ norm as $\Vert g\Vert_2=\sqrt{P g^2}$ and the map $\Omega:L^2(P)\rightarrow L^2(P)$ as $\Omega(g)=\frac{g}{\Vert g\Vert_2}$ if $g\neq 0$. The maximum of the log-likelihood will be denoted :
\[
\lambda_n^k=\sup_{\theta\in \Theta_k}\sum_{i=1}^{n}\log f_{\theta}(z_i)-\log f(z_i).
\]
Finally, let us note $e(z):=\frac{1}{\sigma_0^2}\left(y-\left(\beta^0+\sum_{i=1}^{k^0} a^0_i\phi(b^0_i+{w^0_i}^Tx)\right)\right)$.

For what follows, we will assume the properties:
\begin{description}
\item{H-1 : }the parameters of $\Theta_k$ realizing the true regression function $F_{\theta_0}$ lie in the interior of  $\Theta_k$.
\item{H-2 : }Let $k$ be an integer greater or equal to $k^0$ and 
\[
\theta=\left(\beta,a_1,\cdots,a_k,w_{10},\cdots,w_{1d},\cdots,w_{k0},\cdots,w_{kd}\right).
\]
The model is identifiable in the weak following sense: 
\[
F_{\theta^0}=F_{\theta}\Leftrightarrow\beta^0=\beta\mbox{ and }\sum_{i=1}^{k^0}a^0_i\delta_{w^0_i}=\sum_{i=1}^{k}a_i\delta_{w_i}.
\]
Note that, it is possible that some new constraint on the parameters have to be set to fulfill this assumption. For example, if the transfert function is the hyperbolic tangent (or any odd function), the constraints on the parameters $a_i$ will be : \(a_i\geq\eta\), in order to avoid a symetry on the sign (because $\tanh(-t)=-\tanh(t)$).  
\item{H-3 :} $E(\Vert X\Vert^6)<\infty$.
\item{H-4 :} the functions of the set
\[
\begin{array}{l}
\left(1,\left(x_kx_l\phi^{''}({w^0_i}^Tx)\right)_{1\leq l \leq k\leq d,\ 1\leq i\leq k^0},
\phi^{''}({w^0_i}^Tx)_{1\leq i\leq k^0}, \right.\\
\left. \left(x_k\phi^{'}({w^0_i}^Tx)\right)_{1\leq k\leq d,\ 1\leq i\leq k^0}
 \left(\phi^{'}({w^0_i}^Tx)\right)_{1\leq i\leq k^0},\left(\phi({w^0_i}^Tx)\right)_{1\leq i\leq k^0}\right)
\end{array}
\]
are linearly independent in the Hilbert space $L^2(q\lambda_{d})$. 

\end{description}
We get then the following result:
 
\begin{theorem}
Under the assumptions H-1, H-2 and H-3, a centered Gaussian process $\{W_{S},S\in{\mathbb F}^k\}$ with continuous sample path and covariance kernel $P\left(W_{S_1}W_{S_2}\right)=P\left(S_1S_2\right)$ exists so that 
\[
\lim_{n\rightarrow\infty}2\lambda_n^k=\sup_{S\in{\mathbb F}^k}\left(\max(W_S, 0)\right)^2.
\]
The index set ${\mathbb F}^k$ is defined as ${\mathbb F}^k=\cup_t{\mathbb F}^k_t$, the union runs over any possible $t=\left(t_0,\cdots,t_{k^0}\right)\in{\mathbb N}^{k^0+1}$ with  $0=t_0<t_1<\cdots <t_{k^0}\leq k$ and
\[\begin{array}{l}
{\mathbb F}_t^k=\left\{\Omega\left(\gamma e(z)+\sum_{i=0}^{k^0}\epsilon_ie(z)\phi({w^0_i}^Tx)+\sum_{i=0}^{k^0}e(z)\phi^{'}({w^0_i}^Tx){\zeta}_{i}^Tx\right.\right.\\
\left.+\sum_{i=1}^{k^0}e(z)sg(a^0_i)\phi^{''}({w^0_i}^Tx)\left(\delta(i)\sum_{j=t_{i-1}+1}^{t_i}{\nu_j^t}^Txx^T\nu_j^t\right)\right),\\
\gamma,\epsilon_1,\cdots,\epsilon_{k^0}\in\mathbb R\ ;\ \zeta_1,\cdots,\zeta_{k^0},\left.\nu^t_1,\cdots,\nu^t_{t_{k^0}}\in{\mathbb R}^{d+1}\right\},
\end{array}
\]
where $\delta(i)=1$ if a vector $\bold q$ exists  so that $\sum_{j=t_{i-1}+1}^{t_i}q_j=1$ and $\sum_{j=t_{i-1}+1}^{t_i}\sqrt{q_j}\nu_j^t=0$, otherwise $\delta(i)=0$. The function $sg$ is defined by $sg(x)=1$ if $x>0$ and $sg(x)=-1$ if $x<0$.
\end{theorem}
This theorem is proved in the appendix.
Note that this theorem prove that the LR statistic is tight so penalized likelihood yields a consistent method to identify the minimal architecture of the true model
\section{Identification of the architecture of the MLP}
The point is to guess the number of hidden units of the true MLP function $k^0$. If $k^0$ is known, the information matrix will be regular (see Fukumizu (3)) and pruning of useless parameters will be easy with classical statistical method as in Cottrell et al (2). Here, we assume that the possible number of hidden units in the MLP function is bounded by a large number $K$. So the set of possible parameters will be $\Theta=\cup_{k=1}^{K}\Theta_k$. 

Note that the log-likelihood of the model: $l_n(\theta):=\sum_{i=1}^n\log(f_{\theta}(z_i))$ is known up to the constant $\sum_{i=1}^n\log(x_i)$, independent of the parameter $\theta$. We define $\hat k$, the estimator of maximum of penalized likelihood, as the number of hidden unit maximizing:
\begin{equation}\label{bic}
T_n(k):=\max\{l_n(\theta) : \theta\in\Theta_k\}-p_n(k)
\end{equation}
where $p_n(k)$ is a term which penalizes the log-likelihood in function of the number of hidden units of the model. 

Let $p_n(.)$ be a increasing sequence so that  $p_n(k_1)-p_n(k_2)\rightarrow \infty$ for all $k_1>k_2$ and $\lim_{n\rightarrow\infty}\frac{p_n(k)}{n}=0$. Note that such conditions are verified by BIC-like criterion. 

We get then the following result:
\begin{theorem}
If the assumptions H-1, H-2, H-3 and H-4 are true then  $\hat k\stackrel{P}{\rightarrow}k^0$. 
\end{theorem} 
 The proof is an adaptation of the proof of theorem 2.1 of Gassiat [6]. 
Let us write $ln(f)$ the log-likelihood of the true MLP model.
For any  $f_\theta, \theta\in\Theta$, let
\[
s_\theta(z):=\frac{\frac{f_\theta}{f}(z)-1}{\Vert\frac{f_\theta}{f}-1\Vert_2} 
\mbox{, where } \Vert .\Vert_2 \mbox{ is the } L^2\left(f\lambda_{d+1}\right)\mbox{ norm},
\]
be the generalized score function.
Then, it is obvious to see that under assumptions H-1, H-2, H-3 and H-4 the conditions A1 and A2 of Gassiat (6) are fullfilled. Hence we get the inequality 1.2 of Gassiat (6):
\[
\sup_{\theta\in\Theta}(ln(\theta)-ln(f))\leq\frac{1}{2}\sup_{\theta\in\Theta}\frac{\left(\sum_{i=1}^ns_\theta(z_i)\right)^2}{\sum_{i=1}^n(s_\theta)_{-}^2(z_i)}
\]
where $(s_\theta)_{-}(z)=-\min\left\{0,s_\theta(z)\right\}$.

Now,
\[
\begin{array}{lll}
{\mathbb P}(\hat k > k^0)&\leq&\sum_{k^0+1}^K{\mathbb P}\left(T_n(k)\geq T_n(k^0)\right)\\
&=&\sum_{k^0+1}^K{\mathbb P}\left(\sup_{\theta\in\Theta}(ln(\theta)-ln(f))-\sup_{\theta\in\Theta_{k^0}}(ln(\theta)-ln(f))\right)\\
&\leq&{\mathbb P}\left(\sup_{\theta\in\Theta}\frac{\left(\sum_{i=1}^ns_\theta(Z_i)\right)^2}{\sum_{i=1}^n(s_\theta)_{-}^2(Z_i)}\geq p_n(k)-p_n(k^0)\right)
\end{array}
\]
Now, by Gassiat (6):
\[
\sup_{\theta\in\Theta}\frac{\left(\sum_{i=1}^ns_\theta(Z_i)\right)^2}{\sum_{i=1}^n(s_\theta)_{-}^2(Z_i)}=0_{\mathbb P}(1)
\]
where $0_{\mathbb P}(1)$ means bounded in probability, and
\[
{\mathbb P}(\hat k > k^0)\stackrel{n\rightarrow\infty}{\longrightarrow}0
\]
In the same way
\[
{\mathbb P}(\hat k < k^0)\leq\sum_{k=1}^{k^0-1}{\mathbb P}\left(\sup_{\theta\in\Theta}\frac{(ln(\theta)-ln(f))}{n}\geq\frac{p_n(k)-p_n(k^0)}{n}\right)
\]
But the set $\left\{\log(\frac{f_\theta}{f},\theta\in\Theta\right\}$ is Glivenko-Cantelli, so that $\sup_{\theta\in\Theta}\frac{(ln(\theta)-ln(f))}{n}$ converges in probability to 
\[
-\inf_{\theta\in\Theta}\int\log\frac{f}{f_\theta}<0.
\]
Finally 
\[
{\mathbb P}(\hat k < k^0)\stackrel{n\rightarrow\infty}{\longrightarrow}0
\]
$\blacksquare$
\section{Application to sigmoidal transfert functions}
In this section, The assumptions H-2 and H-4 will be verified for sigmoidal transfert functions :
\[
\phi(t)=\frac{1}{1+e^-t}.
\] 
The assumption H-2 have been shown for hyperbolic tangent functions by Sussmann (8) with additional constraint : \(a_i\geq\eta\), morevover MLP with sigmoidale tranfert functions or hyperbolic tangente transfert functions are equivalent, because an one-to-on correspondence between the two kinds of MLP exists as \(\frac{1}{1+e^-t}=(1+\tanh(t/2))/2\) (see Fukumizu (3)).Hence the assumption H-2 is verified for sigmoidale functions with the additional constraint. 

The main point is to verify H-4. The proof use an extension of the result of Fukumizu[3]. 

We define the complex sigmoidal function on $\mathbb C$ by $\phi(z)=\frac{1}{1+e^{-z}}$.

The singularities of $\phi$ are : 
\[
\left\{z\in\mathbb C\left|z=(2n+1)\pi\sqrt{-1},n\in\mathbb Z\right.\right\}
\]
all of which are poles of order 1.
Next we review a fundamental propositions in complex analysis.
\begin{e-proposition}
Let $\phi$ be a holomorphic function on a connected open set $D$ in $\mathbb C$ and $p$ be a point in $D$. If a sequence $\left\{p_n\right\}_{n=1}^{\infty}$ exists in $D$ so that $p_n\neq p,\lim_{n\rightarrow n}p_n=p$ and $\phi(p_n)=0$ for all $n\in\mathbb N$ then $\phi(z)=0$ for all $z\in D$.
\end{e-proposition}
\begin{e-proposition}\label{holo}
Le $\phi$ be a holomorphic function on a connected open set $D$ in $\mathbb C$, and $p$ be a point in $D$. 
Then the following equivalence relations hold: 
\begin{itemize}
\item $p$ is a removable singularity
\[
\Leftrightarrow \lim_{z\rightarrow p}f(z)\in \mathbb C
\]
\item $p$ is a pole
\[
\Leftrightarrow \lim_{z\rightarrow p}\left|f(z)\right|=\infty
\]
\item $p$ is an essential singularity
\[
\Leftrightarrow \lim_{z\rightarrow p}\left|f(z)\right|\mbox{ does not exist}
\]
\end{itemize}
\end{e-proposition}
Let $w_0,\cdots,w_{k^0}$ be the parameters of the minimale, true, MLP function (in order to simplify the notations,  the exponent ``0'' is missing). By the lemma 3 of Fukumizu (3), a basis of ${\mathbb R}^d$ $\left(x^{(1)},\cdots,x^{(d)}\right)$ exists so that
\begin{enumerate}
\item For all $i\in \{1,\cdots,k^0\}$ and all $h\in \{1,\cdots,d\}$ 
\[\sum_{j=1}^d{w_{ij}}x^{(h)}_j\neq 0.\]
\item For all $i_1,i_2\in \{1,\cdots,k^0\}$, $i_1\neq i_2$ and all $h\in \{1,\cdots,d\}$ 
\[
w_{i_10}+\sum_{j=1}^dw_{i_1j}x^{(h)}_j\neq \pm \left(w_{i_20}+\sum_{j=1}^dw_{i_2j}x^{(h)}_j\right)
\]
\end{enumerate}
 For $h$, $1\leq h\leq d$ and $i\in \{1,\cdots,k^0\}$ let be $m^{(h)}_i:=\sum_{j=1}^d{w_{ij}}x^{(h)}_j$. 
We fix $l$ for a while. We set
\[
S^{(l)}_i=\left\{ u\in{\mathbb C}\left|u=\frac{(2n+1)\pi \sqrt{-1}-w_{i0}}{m^{(l)}_i}, n\in {\mathbb Z}\right.\right\}
\]
Clearly the points in $S^{(l)}_i$ are the singularities of $\phi\left(m_i^{(l)}u+w_{i0}\right)$.
Note that these points are pole of order 1 for 
\[\phi(m_i^{(l)}u+w_{i0})=\frac{1}{1+e^{-\left(m_i^{(l)}u+w_{i0}\right)}}
\]
of order 2 for 
\[\phi'(m_i^{(l)}u+w_{i0})=-\frac{e^{-\left(m_i^{(l)}u+w_{i0}\right)}}{\left(1+e^{-\left(m_i^{(l)}u+w_{i0}\right)}\right)^2}
\] 
and 3 for 
\[\phi''(m_i^{(l)}u+w_{i0})=\frac{e^{-\left(m_i^{(l)}u+w_{i0}\right)}}{\left(1+e^{-\left(m_i^{(l)}u+w_{i0}\right)}\right)^2}+2\frac{e^{-2\left(m_i^{(l)}u+w_{i0}\right)}}{\left(1+e^{-\left(m_i^{(l)}u+w_{i0}\right)}\right)^3}
\]
Let be $D^(l):={\mathbb C}-\cup_{1\leq i\leq k^0}S_i^{(l)}$, Holomorphic functions on $D^{(l)}$ are defined as follows:
\[
\begin{array}{l}
\Psi^{(l)}(u):=\alpha_0+\sum_{i=1}^{k^0}\alpha_i\phi(m_i^{(l)}u+w_{i0})+\sum_{i=1}^{k^0}\epsilon_i\phi^{'}(m_i^{(l)}u+w_{i0})\\
+\sum_{i=1}^{k^0}\sum_{j=1}^d\beta_{ij}\phi^{'}(m_i^{(l)}u+w_{i0})x_j^{(l)}u
+\sum_{i=1}^{k^0}\delta_i\phi^{''}(m_i^{(l)}u+w_{i0})\\
+\sum_{i=1}^{k^0}\sum_{j,k=1,\ j\leq k}^d\gamma_{ijk}\phi^{''}(m_i^{(l)}u+w_{i0})x_j^{(l)}x_k^{(l)}u^2
\end{array}
\]
The functions in the set
\[
\begin{array}{l}
\left(1,\left(x_kx_l\phi^{''}({w^0_i}^Tx)\right)_{1\leq l \leq k\leq d,\ 1\leq i\leq k^0},\phi^{''}({w^0_i}^Tx)_{1\leq i\leq k^0}, \left(x_k\phi^{'}({w^0_i}^Tx)\right)_{1\leq k\leq d,\ 1\leq i\leq k^0}\right.\\
\left.  \left(\phi^{'}({w^0_i}^Tx)\right)_{1\leq i\leq k^0},\left(\phi({w^0_i}^Tx)\right)_{1\leq i\leq k^0},\left(\phi({w_i}^Tx)\right)_{1\leq i\leq l}\right)
\end{array}
\]
are linearly independent if the following property is verified :
\[
 \forall u\in D(l), \Psi^{(l)}(u)=0 \Leftrightarrow \alpha_i,\epsilon_i,\beta_{ij},\delta_i\mbox{ and }\gamma_{ijk}\mbox{ are equal to 0}
\]
Let us assume that: $\forall u\in D^{(l)}, \Psi^{(l)}(u)=0$, then by proposition \ref{holo} all the point in $S^{(l)}_i$ are removable singularities.

Let us write
\[
p_i^{(l)}:=\frac{\pi\sqrt{-1}-b^0_i}{m_i^{(l)}}\in S_i^{(l)}
\]

Clearly , for $1\leq i\leq k^0-1$, $p_{k^0}^{(l)}\notin S_i^{(l)}$, because for all $i_1,i_2\in \{1,\cdots,k^0\}$, $i_1\neq i_2$ and all $h\in \{1,\cdots,d\}$ 
\[
w_{i_10}+\sum_{j=1}^dw_{i_1j}x^{(h)}_j\neq \pm \left(w_{i_20}+\sum_{j=1}^dw_{i_2j}x^{(h)}_j\right)
\]
So,  $\Psi^{(l)}(u)$ can be written as:
\[
\begin{array}{l}
\Psi^{(l)}(u)=\alpha_{k^0}\phi(m_i^{(l)}u+w_{i0})+\left(\sum_{i=1}^{d}\beta_{k^0i}x_i^{(l)}u+\epsilon_{k^0}\right)\phi^{'}(m_{k^0}^{(l)}u+w_{k^00})\\
+\left(\sum_{i,j=1,\ i\leq j}^{d}\gamma_{k^0ij}x_i^{(l)}x_j^{(l)}u^2+\delta_1\right)\phi^{''}(m_{k^0}^{(l)}u+w_{k^00})+\Psi^{(l)}_{k^0-1}(u)
\end{array}
\]
where
\[
\begin{array}{l}
\Psi^{(l)}_{k^0-1}(u):=\alpha_0+\sum_{i=1}^{k^0-1}\alpha_i\phi(m_i^{(l)}u+w_{i0})\\
+\sum_{i=1}^{k^0-1}\epsilon_i\phi^{'}(m_i^{(l)}u+w_{i0})+\sum_{i=1}^{k^0-1}\sum_{j=1}^d\beta_{ij}\phi^{'}(m_i^{(l)}u+w_{i0})x_j^{(l)}u\\
+\sum_{i=1}^{k^0-1}\delta_i\phi^{''}(m_i^{(l)}u+w_{i0})+\sum_{i=1}^{k^0-1}\sum_{j,k=1}^d\gamma_{ijk}\phi^{''}(m_i^{(l)}u+w_{i0})x_j^{(l)}x_k^{(l)}u^2
\end{array}
\]
The point $p_{k^0}^{(l)}$ is a regular point of $\Psi^{(l)}_{{k^0}-1}(u)$ while $\phi(m_{k^0}^{(l)}u+w_{{k^0}0})$ has a pole of order 1 at $p_{k^0}^{(l)}$, $\phi^{'}(m_{k^0}^{(l)}u+w_{k^00})$ has a pole of order 2 at $p_{k^0}^{(l)}$ and $\phi^{''}(m_{k^0}^{(l)}u+w_{k^00})$ has a pole of order 3 et $p_{k^0}^{(l)}$. Since $p_{k^0}^{(l)}$ is a removable singularity of $\Psi^{(l)}(u)$, we have:
\[
\alpha_{k^0}=0,\ \epsilon_{k^0}=0,\ \sum_{i=1}^{d}\beta_{k^0i}x_i^{(l)}=0\mbox{ and }\sum_{i,j=1,\ i\leq j}^{d}\gamma_{k^0ij}x_i^{(l)}x_j^{(l)}=\delta_{k^0}=0
\]
As a result  $\Psi^{(l)}(u)=\Psi^{(l)}_{k^0-1}(u)$. Applying the same argument successively to $p_{k^0-1}^{(l)},\cdots,p_1^{(l)}$ we finally obtain, for all $1\leq i\leq k^0$, $1\leq j\leq k\leq d$:
\[
\begin{array}{c}
\alpha_i=0\\
\epsilon_i=0\\
\sum_{j=1}^{d}\beta_{ij}x_j^{(l)}=0\\
\sum_{j,k=1,\ j\leq k}^{d}\gamma_{ijk}x_j^{(l)}x_k^{(l)}=0\\
\delta_i=0\\
\alpha_0=0
\end{array}
\]
Since $\left(x^{(1)},\cdots,x^{(d)}\right)$ form a basis of ${\mathbb R}^d$, we have  $\beta_{ij}=0$ for all $1\leq i\leq k^0$, $1\leq j \leq d$. 

For $\gamma_{ijk}$, we get:
\[
\sum_{j,k=1,\ j\leq k}^{d}\gamma_{ijk}x_j^{(l)}x_k^{(l)}=\sum_{k=1}^d\left(\sum_{j=1}^k\gamma_{ijk}x_j^{(l)}\right)x_k^{(l)}=0
\]
and, since $\left(x^{(1)},\cdots,x^{(d)}\right)$ form a basis of ${\mathbb R}^d$, we obtain, for all $l\in\{1,\cdots,d\}$: 
\[
\begin{array}{c}
\gamma_{i11}x_1^{(l)}=0\\
\vdots\\
\sum_{j=1}^k\gamma_{ijk}x_j^{(l)}=0\\
\vdots\\
\sum_{j=1}^d\gamma_{ijd}x_j^{(l)}=0\\
\end{array}
\]
and by the same remark on $\left(x^{(1)},\cdots,x^{(d)}\right)$,  $\gamma_{ijk}=0$ for all $1\leq i\leq k^0$, $1\leq j\leq k \leq d$.

This prove that H-4 hold for sigmoidale functions with the additional constraints \(\forall i\in\{1,\cdots,k^0\}, a_i\geq\eta\)$\blacksquare$   

\section{Conclusion}
We have computed the asymptotic distribution of the LR statistic for parametric MLP regression. 
This theorem can be applied to the most widely used transfer functions for MLP: the sigmoidal functions.Note that the results assume some constraints on the parameters of the MLP, the constraints on $a_i$ and $w_i$ may be relaxed, but a more clever reparameterization and a higher order in the development of the LR statistics should certainly be required. The asumption on $\sigma^2$ is certainly easier to remove, however the development of lemma \ref{dev} will be much more complicated and so the limit score functions. Finally, this theorem shows that the LR statistic is tight, so information criteria such as the Bayesian information criteria (BIC) will be consistent in the sense that they will select the model with the true dimension $k^0$ with probability 1, as the number of observations goes to infinite. This is the main pratical application of the results obtained in the paper.

\section{Appendix: Proof of the Theorem.}
Let
\[
s_\theta(z):=\frac{\frac{f_\theta}{f}(z)-1}{\Vert\frac{f_\theta}{f}-1\Vert_2} 
\mbox{, where } \Vert .\Vert_2 \mbox{ is the } L^2\left(f\lambda_{d+1}\right)\mbox{ norm},
\]
be the generalized score functions. Firstly, we will get an asymptotic development of the generalized score when the model is over-parameterized. We will reparameterize the model using the same method as Liu and Shao [7] for the mixture models.
\subsection{Reparametrization.} 

If $\frac{f_\theta}{f}-1=0$, we have $\beta=\beta^0$ and a vector $t=(t_i)_{1\leq i\leq k^0}$ exists so that $0=t_0<t_1<\cdots<t_{k^0}\leq k$ and, up to permutations, we have \(w_{t_{i-1}+1}=\cdots=w_{t_i}=w^0_i\), \(\sum_{j=t_{i-1}+1}^{t_i}a_j=a_i^0\). Let $s_i=\sum_{j=t_{i-1}+1}^{t_i}a_j-a_i^0$ be and $q_j=\frac{a_j}{\sum_{t_{i-1}+1}^{t_i}a_j}$ if $\sum_{t_{i-1}+1}^{t_i}a_j\neq 0$ and otherwise $q_j=0$, we get then the reparameterization $\theta=\left(\Phi_t,\psi_t\right)$ with 
\[
\Phi_t=\left(\beta,(w_j)_{j=1}^{t_{k^0}},(s_i)_{i=1}^{k^0}\right),\ \psi_t=\left((q_j)_{j=1}^{t_{k^0}}\right).
\]
With this parameterization, for a fixed $t$,  $\Phi_t$ is an identifiable parameter and all the non-identifiability of the model will be in $\psi_t$. Then $\frac{f_\theta}{f}(z)$ will be equal to
\[
\frac{exp\left(-\frac{1}{2\sigma^2}\left(y-\left(\beta+\sum_{i=1}^{k^0}(s_i+a^0_i)\sum_{j=t_{i-1}+1}^{t_i}q_j\phi(w_j^Tx)\right)\right)^2\right)}{exp\left(-\frac{1}{2\sigma_0^2}\left(y-\left(\beta^0+\sum_{i=1}^{k^0} a^0_i\phi({w^0_i}^Tx)\right)\right)^2\right)}.
\]
Now, as the third derivative of the transfer function is bounded and thanks to the assumption H-2, the third order derivative of the function $\frac{f_\theta}{f}(z)$ with respect to the components of $\Phi_t$ will be dominated by a square integrable function, because there exists a constant $C$ so that we have the following inequalities:
\[
\forall\theta_i,\theta_j,\theta_l\in\left\{w_{10},\cdots,w_{kd}\right\},\ \sup_{\theta\in\Theta_k}\Vert \frac{\partial^3 F_\theta(X)}{\partial \theta_i\partial \theta_j\partial \theta_l}\Vert\leq C(1+\Vert X\Vert^3).
\]
So, by the Taylor formula with an integral remainder around the identifiable parameter $\Phi^0_t$ with

\[
\begin{array}{cccccc}
\Phi^0_t=(\beta^0,&\underbrace{w_1^0,\cdots,w_1^0}&,\cdots,&\underbrace{w_{k^0}^0,\cdots,w_{k^0}^0}&,\underbrace{0,\cdots,0}),\\
&t_1& &t_{k^0}-t_{k^0-1}&k^0
\end{array}
\]
we get the following Taylor expansion for the likelihood ratio :    
\begin{lemma}\label{dev}
For a fixed $t$, let us write $D(\Phi_t,\psi_t):=\Vert\frac{f_{(\Phi_t,\psi_t)}}{f}-1\Vert_2$. In the neighborhood  of the identifiable parameter $\Phi^0_t$, we get the following approximation:
\[
\frac{f_\theta}{f}(z)=1+(\Phi_t-\Phi^0_t)^Tf^{'}_{(\Phi^0_t,\psi_t)}(z)+0.5(\Phi_t-\Phi^0_t)^Tf^{''}_{(\Phi^0_t,\psi_t)}(z)(\Phi_t-\Phi^0_t)+o(D(\Phi_t,\psi_t)),
\]
with
\[
\begin{array}{l}
(\Phi_t-\Phi^0_t)^Tf^{'}_{(\Phi^0_t,\psi_t)}(z)=
e(z)\left(\beta-\beta^0+\sum_{i=1}^{k^0}s_i\phi({w^0_i}^Tx)\right.\\
\left.+\sum_{i=1}^{k^0}\sum_{j=t_{i-1}+1}^{t_i}q_j\left(w_{j}-w^0_{i}\right)^Txa^0_i\phi^{'}({w^0_i}^Tx)\right)
\end{array}
\]
and
\[
\begin{array}{l}
(\Phi_t-\Phi^0_t)^Tf^{''}_{(\Phi^0_t,\psi_t)}(z)(\Phi_t-\Phi^0_t)=\\
\left(1-\frac{1}{e^2(z)}\right)\left((\Phi_t-\Phi^0_t)^Tf^{'}_{(\Phi^0_t,\psi_t)}(z){f^{'}_{(\Phi^0_t,\psi_t)}}^T(z)(\Phi_t-\Phi^0_t)\right)\\
+e(z)\times\left(\sum_{i=1}^{k^0}\sum_{j=t_{i-1}+1}^{t_i}q_j(w_{j}-w^0_{i})^Txx^T(w_{j}-w^0_{i})a^0_i\phi^{''}({w^0_i}^Tx)\right.\\
\left. +\sum_{i=1}^{k^0}\sum_{j=t_{i-1}+1}^{t_i}(q_jw_{j}-w^0_{i})^Txs_{i}\phi^{'}({w^0_i}^Tx) \right).
\end{array}
\]
\end{lemma}
\paragraph{Proof of the lemma.}

This development, obtained by a straightforward calculation of the derivatives of $\frac{f_\theta}{f}(z)$ with respect to the components of $\Phi_t$ up to the second order, is postponed to the end of this appendix.
 
Now, the convergence to a Gaussian process will be derived from the Donsker property of the set of generalized score functions $\mathbb S=\left\{s_\theta(z),\theta\in\Theta_k\right\}$.Let an $\varepsilon$-bracket $[l,u]$ be a set of function $h$ with $l\leq h\leq u$ with $\sqrt{P(l-u)^2}<\varepsilon$. The bracketing number $N_{[]}\left(\varepsilon,\mathbb S,L^2\left(f\lambda_{d+1}\right)\right)$ is the minimum number of $\varepsilon$-brackets needed to cover $\mathbb S$. The entropy with bracketing is the logarithm of the bracketing number. It is well known (see van der Vaart [8]) that the class of functions $\mathbb S$ will be Donsker if its entropy with bracketing grows with a slower order than $\frac{1}{\varepsilon}^2$. A sufficient condition for Donsker property is then that the bracketing number grows as a polynomial function of $\frac{1}{\varepsilon}$. 
\subsection{Polynomial bound for the growth of  bracketing number.}

Let us write $D(\theta):=\Vert\frac{f_{(\theta)}}{f}-1\Vert_2$, for all $\varepsilon>0$, the set of parameters can be divided in two sets: ${\mathbb S}_\varepsilon$ and ${\mathbb S}_0$ with
\[
  {\mathbb S}_\varepsilon=\left\{\theta\in \Theta_k\mbox{ so that }D(\theta)\geq\varepsilon\right\}\mbox{ and }{\mathbb S}_0=\left\{\theta\in \Theta_k\mbox{ so that }D(\theta)<\varepsilon\right\}.
\] 

For $\theta_1$ and $\theta_2$ belonging to ${\mathbb S}_\varepsilon$, we get:
\[
\begin{array}{l}
\left\Vert\frac{\frac{f_{\theta_1}}{f}-1}{\left\Vert\frac{f_{\theta_1}}{f}-1\right\Vert_2}-\frac{\frac{f_{\theta_2}}{f}-1}{\left\Vert\frac{f_{\theta_2}}{f}-1\right\Vert_2}\right\Vert_2=\left\Vert\frac{\frac{f_{\theta_1}}{f}-1}{\left\Vert\frac{f_{\theta_1}}{f}-1\right\Vert_2}-\frac{\frac{f_{\theta_2}}{f}-1}{\left\Vert\frac{f_{\theta_1}}{f}-1\right\Vert_2}+\frac{\frac{f_{\theta_2}}{f}-1}{\left\Vert\frac{f_{\theta_1}}{f}-1\right\Vert_2}-\frac{\frac{f_{\theta_2}}{f}-1}{\left\Vert\frac{f_{\theta_2}}{f}-1\right\Vert_2}\right\Vert_2\\
\leq\left\Vert\frac{\frac{f_{\theta_1}}{f}-1}{\left\Vert\frac{f_{\theta_1}}{f}-1\right\Vert_2}-\frac{\frac{f_{\theta_2}}{f}-1}{\left\Vert\frac{f_{\theta_1}}{f}-1\right\Vert_2}\right\Vert_2+\left\Vert\frac{\frac{f_{\theta_2}}{f}-1}{\left\Vert\frac{f_{\theta_1}}{f}-1\right\Vert_2}-\frac{\frac{f_{\theta_2}}{f}-1}{\left\Vert\frac{f_{\theta_2}}{f}-1\right\Vert_2}\right\Vert_2\leq2\frac{\left\Vert\frac{f_{\theta_1}}{f}-\frac{f_{\theta_2}}{f}\right\Vert_2}{\left\Vert\frac{f_{\theta_1}}{f}-1\right\Vert_2}\leq2\frac{\left\Vert\frac{f_{\theta_1}}{f}-\frac{f_{\theta_2}}{f}\right\Vert_2}{\varepsilon}.
\end{array}
\]
Hence, on ${\mathbb S}_\varepsilon$, it is sufficient that
\[
\left\Vert\frac{f_{\theta_1}}{f}-\frac{f_{\theta_2}}{f}\right\Vert_2<\frac{\varepsilon^2}{2}
\]
for
\[
\left\Vert\frac{\frac{f_{\theta_1}}{f}-1}{\left\Vert\frac{f_{\theta_1}}{f}-1\right\Vert_2}-\frac{\frac{f_{\theta_2}}{f}-1}{\left\Vert\frac{f_{\theta_2}}{f}-1\right\Vert_2}\right\Vert_2<\varepsilon.
\]
Now, ${\mathbb S}_\varepsilon$ is a parametric class. Since the derivatives of the transfer functions are bounded and $E\Vert X\Vert<\infty$ a function $m(z)$ exists , with $E[m(z)]<\infty$, so that  
\[
\forall\theta_i\in \left\{\beta,a_1,\cdots,a_k,w_{10},\cdots,w_{1d},\cdots,w_{kd}\right\},\ 
\left|\frac{\partial \frac{f_\theta}{f}}{\partial \theta_i}(z)\right|\leq m(z). 
\]
According to the exemple 19.7 of van der Vaart [8], it exists a constant $K$ so that the bracketing number of ${\mathbb S}_\varepsilon$ is lower than 
\[
K\left(\frac{\mbox{diam}\Theta_k}{\varepsilon^2}\right)^{k\times(d+2)+1}=K\left(\frac{\sqrt{\mbox{diam}\Theta_k}}{\varepsilon}\right)^{k\times(2d+4)+2},
\]
where $\mbox{diam}\Theta_k$ is the diameter of the smallest sphere of $\mathbb R^k$ including $\Theta_k$. 

For $\theta$ belonging to ${\mathbb S}_0$, \(\frac{f_\theta(z)}{f}-1\) is the sum of a linear combination of
\[
\begin{array}{l}
V(z):=\left(e(z),\left(e(z)x_kx_l\phi^{''}({w^0_i}^Tx)\right)_{1\leq l \leq k\leq d,\ 1\leq i\leq k^0},e(z)\phi^{''}({w^0_i}^Tx)_{1\leq i\leq k^0},\right.\\
\left(e(z)x_k\phi^{'}({w^0_i}^Tx)\right)_{1\leq k\leq d,\ 1\leq i\leq k^0}, \left. \left(e(z)\phi^{'}({w^0_i}^Tx)\right)_{1\leq i\leq k^0},\left(e(z)\phi({w^0_i}^Tx)\right)_{1\leq i\leq k^0}\right)
\end{array}
\] 

and of a term whose $L^2\left(f\lambda_{d+1}\right)$ norm is negligible compared to the $L^2\left(f\lambda_{d+1}\right)$ norm of this combination when $\varepsilon$ goes to 0. 
By assumption H-3, a strictly positive number $m$ exists so that for any vector of norm 1 with components
\[
C=\left(c,c_1,\cdots,c_{k^0\times\frac{d(d+1)}{2}},d_1,\cdots,d_{k^0},e_1,\cdots,e_{k^0\times d},f_1,\cdots,f_{k^0},g_1,\cdots,g_{k^0}\right)
\]
and $\varepsilon$ sufficiently small:
\[
\Vert C^TV(z)\Vert_2> m+\varepsilon.
\]
Since any function $\frac{\frac{f_\theta}{f}-1}{\Vert \frac{f_\theta}{f}-1\Vert_2}$ can be written:
\[
\frac{C^TV(z)+o(\Vert C^TV(z)\Vert_2)}{\Vert C^TV(z)+o(\Vert C^TV(z)\Vert_2)\Vert_2},
\] 
$\mathbb S_0$ belongs to the set of functions:
\[
\left\{D^TV(z)+o(1),\Vert D\Vert_2\leq\frac{1}{m}\right\}\subset\left\{D^TV(z)+\gamma,\Vert D\Vert_2\leq\frac{1}{m},|\gamma|<1\right\}
\] 
whose bracketing number is smaller or equal to $O\left(\frac{1}{\varepsilon}\right)^{k^0\times\left(\frac{d(d+1)}{2}+d+3\right)+2}$.

This proves that the bracketing number of $\mathbb S$ is polynomial, hence $\mathbb S$ is a Donsker class.
\subsection{Asymptotic index set.}

Since the class of generalized score functions $\mathbb S$ is a Donsker class the theorem follows from theorem 3.1 of Gassiat [6] or theorem 3.1 of Liu and Shao [7]. Following these authors, the set of limit score functions $\mathbb F^k$ is defined as the set of functions  $d$ so that one can find a sequence $g_n:=f_{\theta_k^n}, \theta_k^n\in \Theta_k$ satisfying $\Vert \frac{g_n-f}{f}\Vert_2\rightarrow 0$ and $\Vert d-s_{g_n}\Vert_2\rightarrow 0$, where $s_{g_n}=\frac{\frac{g_n}{f}(z)-1}{\Vert\frac{g_n}{f}-1\Vert_2}$. Note that, for a particular sequence of maximum likelihood estimators $(\theta^n)_{n\in\mathbb{N}}$, the partition of the indices can depend on $n$, but $(\theta^n)_{n\in\mathbb{N}}$ will be the union of converging sub-sequences belonging the set of limit score functions.   

Let us define the two principal behaviors for the sequences $g_n$ which influence the form of functions $d$ :
\begin{itemize}
\item If the second order term is negligible behind the first one :
\[
\frac{f_{\theta_k^n}}{f}(z)-1=(\Phi_k^n-\Phi^0)^Tf^{'}_{(\Phi^0_t,\psi_k^n)}(z)+o(D(\Phi_k^n,\psi_k^n)).
\]
\item If the second order term is not negligible compared to the first one :
\[
\begin{array}{l}
\frac{f_{\theta_k^n}}{f}(z)-1=(\Phi_k^n-\Phi^0)^Tf^{'}_{(\Phi^0_t,\psi_k^n)}(z)+\\
0.5(\Phi_k^n-\Phi^0_t)^Tf^{''}_{(\Phi^0_t,\psi_k^n)}(z)(\Phi_k^n-\Phi^0_t)+o(D(\Phi_k^n,\psi_k^n)).
\end{array}
\]
\end{itemize}
In the first case, each sequence $g_n$ is the finite union of convergent subsequences $g_k(n)$ and for each subsequence  a set $t=\left(t_0,\cdots,t_{k^0}\right)\in{\mathbb N}^{k^0+1}$ (with  $0=t_0<t_1<\cdots <t_{k^0}\leq k$) exists so that the limit functions $d$ of  $s_{g_k(n)}$ will be:
\[\begin{array}{ll}
\mathbb D^t_1=&\Omega\left\{\gamma e(z)+\sum_{i=0}^{k^0}\epsilon_ie(z)\phi({w^0_i}^Tx)+\sum_{i=0}^{k^0}e(z)\phi^{'}({w^0_i}^Tx){\zeta}_{i}^Tx\right.\\
&\left.\gamma,\epsilon_1,\cdots,\epsilon_{k^0}\in\mathbb R\ ;\ \zeta_1,\cdots,\zeta_{k^0}\in{\mathbb R}^{d+1}\right\}.
\end{array}
\]

In the second case, each sequence $g_n$ is the finite union of convergent subsequences $g_k(n)$ and for each subsequence, an index $i$ exists so that :
\[
\sum_{j=t_{i-1}+1}^{t_i}q_j(w_{j}-w^0_{i})=0,
\]
otherwise the second order term will be negligible compared to the first one, so
\[
\sum_{j=t_{i-1}+1}^{t_i}\sqrt{q_j}\times\sqrt{q_j}(w_{j}-w^0_{i})=0.
\]

Hence, a set $t=\left(t_0,\cdots,t_{k^0}\right)\in{\mathbb N}^{k^0+1}$ exists, with  $0=t_0<t_1<\cdots <t_{k^0}\leq k$ so that the set of functions $d$ will be: 
\[\begin{array}{l}
\Omega\left\{\gamma e(z)+\sum_{i=0}^{k^0}\epsilon_ie(z)\phi({w^0_i}^Tx)+\sum_{i=0}^{k^0}e(z)\phi^{'}({w^0_i}^Tx){\zeta}_{i}^Tx\right.\\
+\sum_{j=t_{k^0+1}}^k\alpha^t_je(z)\phi({w_j}^Tx)\\
+\sum_{i=1}^{k^0}e(z)sg(a^0_i)\phi^{''}({w^0_i}^Tx)\left(\delta(i)\sum_{j=t_{i-1}+1}^{t_i}{\nu_j^t}^Txx^T\nu_j^t\right):\\
\gamma,\epsilon_1,\cdots,\epsilon_{k^0},\alpha^t_{k^0+1},\cdots,\alpha^t_{k}\in\mathbb R\ ;\ \zeta_1,\cdots,\zeta_{k^0}\in{\mathbb R}^{d+1}\\
\left.\mu^t_1,\cdots,\mu^t_{t_{k^0}}\nu^t_1,\cdots,\nu^t_{t_{k^0}}\in{\mathbb R}^d\right\}\supset\mathbb D^t_1,
\end{array}
\]
where $\delta(i)=1$ if a vector $\bold q$ exists with $\sum_{j=t_{i-1}+1}^{t_i}q_j=1$ and  $\sum_{j=t_{i-1}+1}^{t_i}\sqrt{q_j}\nu_j=0$, otherwise $\delta(i)=0$.

So, the limit functions $d$ will belong to $\mathbb F^k$. 

Conversely, for  $x\in L^2(\lambda_{d+1})$, let $d$ be an element of $\mathbb F^k$:
\[\begin{array}{l}
d=\Omega\left(\gamma e(z)+\sum_{i=0}^{k^0}\epsilon_ie(z)\phi({w^0_i}^Tx)+\sum_{i=0}^{k^0}e(z)\phi^{'}({w^0_i}^Tx){\zeta}_{i}^Tx\right.\\
\left.+\sum_{i=1}^{k^0}e(z)sg(a^0_i)\phi^{''}(b^0_i+{w^0_i}^Tx)\left(\delta(i)\sum_{j=t_{i-1}+1}^{t_i}{\nu_j^t}^Txx^T\nu_j^t\right)\right).
\end{array}
\]
As functions $d$ belong to the  Hilbert sphere, one of their components is not equal to 0. Let us assume that this component is $\gamma$, but the proof would be similar with any other component.  The norm of  $d$ is $1$, so any component of $d$ is determined by the ratio: $\frac{\epsilon_1}{\gamma},\cdots,\frac{1}{\gamma}\nu_{k^0}^t$.

Then, we can chose $\theta_k^n=\left(\beta^n,a^n_1,\cdots,a^n_k,w^n_{10},\cdots,w^n_{1d},\cdots,w^n_{kd}\right)$ so that:
\[
\begin{array}{l}
\forall i\in\{1,\cdots,k^0\}\ :\ \frac{s^n_i}{\beta_n-\beta^0}\stackrel{n\rightarrow \infty}{\longrightarrow}\frac{\epsilon_i}{\gamma},\\
\forall i\in\{1,\cdots,k^0\}\ :\ \sum_{j=t_{i-1}+1}^{t_i}\frac{q^n_j}{\beta_n-\beta^0}\left(w^n_j-w_i^0\right)\stackrel{n\rightarrow \infty}{\longrightarrow}\frac{1}{\gamma}\zeta_i,\\
\forall j\in\{1,\cdots,t_{k^0}\}\ :\ \frac{\sqrt{q^n_j}}{\beta_n-\beta^0}\left(w^n_j-w_i^0\right)\stackrel{n\rightarrow \infty}{\longrightarrow}\frac{1}{\gamma}\nu_j,
\end{array}
\] 
since $\Theta_k$ contains a neighborhood of the parameters realizing the true regression function $F_{\theta^0}$.
\(\blacksquare\)
\subsection{The derivatives of the LR statistic}

\subsubsection{Calculation of $(\Phi_t-\Phi^0_t)^Tf^{'}_{(\Phi^0_t,\psi_t)}(z)$}
In the sequel, we write $x:=\left(1,x_1,\cdots,x_d\right)^T$.

To get \((\Phi_t-\Phi^0_t)^Tf^{'}_{(\Phi^0_t,\psi_t)}(z)\), we compute the derivatives of the $\frac{f_\theta}{f}(z)$ with respect to each parameter of
\(
\Phi_t=\left(\beta,(w_j)_{j=1}^{t_{k^0}},(s_i)_{i=1}^{k^0}\right)
\).
 
Let us recall that $e(z):=\frac{1}{\sigma_0^2}\left(y-\left(\beta^0+\sum_{i=1}^{k^0} a^0_i\phi(b^0_i+{w^0_i}^Tx)\right)\right)$, we get: 
\begin{itemize}

\item \[
\frac{\partial \frac{f_\theta}{f}(z)}{\partial \beta}(\Phi^0_t)=e(z)
\]
\item \[
\frac{\partial \frac{f_\theta}{f}(z)}{\partial s_i}(\Phi^0_t)=e(z)\sum_{j=t_{i-1}+1}^{t_i}q_j\phi({w^0_j}^Tx)=e(z)\phi({w^0_i}^Tx)
\]

\item For $j\in\{t_{i-1}+1,\cdots,t_i\}$, let us write $\frac{\partial \frac{f_\theta}{f}(z)}{\partial w_j}:=\left(\frac{\partial \frac{f_\theta}{f}(z)}{\partial w_{j0}},\cdots,\frac{\partial \frac{f_\theta}{f}(z)}{\partial w_{jd}}\right)^T$\[
\frac{\partial \frac{f_\theta}{f}(z)}{\partial w_j}(\Phi^0_t)=e(z)a^0_iq_j\phi^{'}({w^0_i}^Tx)x
\]

\item For $j\in\{t_{k^0}+1,\cdots,k\}$ :
\[
\frac{\partial \frac{f_\theta}{f}(z)}{\partial a_j}=e(z)\phi({w_j}^Tx)
\]
 
\end{itemize}

These equations yield us the expression of $(\Phi_t-\Phi^0_t)^Tf^{'}_{(\Phi^0_t,\psi_t)}(z)$.

\subsubsection{Calculation of $(\Phi_t-\Phi^0_t)^Tf^{''}_{(\Phi^0_t,\psi_t)}(z)(\Phi_t-\Phi^0_t)$}

\begin{itemize}

\item \[
\frac{\partial^2\frac{f_\theta}{f}(z)}{\partial \beta^2}(\Phi^0_t)=e^2(z)-1
\]
\item \[
\frac{\partial^2\frac{f_\theta}{f}(z)}{\partial \beta\partial s_i}(\Phi^0_t)=\left(e^2(z)-1\right)\phi({w^0_i}^Tx)
\]

\item For $j\in\{t_{i-1}+1,\cdots,t_i\}$, let us write $\frac{\partial^2\frac{f_\theta}{f}(z)}{\partial\beta\partial w_j}:=\left(\frac{\partial^2 \frac{f_\theta}{f}(z)}{\partial\beta\partial w_{j0}},\cdots,\frac{\partial^2\frac{f_\theta}{f}(z)}{\partial\beta\partial w_{jd}}\right)^T$\[
\frac{\partial^2\frac{f_\theta}{f}(z)}{\partial \beta\partial w_j}(\Phi^0_t)=\left(e^2(z)-1\right)a^0_iq_j\phi^{'}({w^0_i}^Tx)x
\]

\item For $j\in\{t_{k^0}+1,\cdots,k\}$ :
\[
\frac{\partial^2\frac{f_\theta}{f}(z)}{\partial \beta\partial a_j}(\Phi^0_t)=\left(e^2(z)-1\right)\phi({w_j}^Tx)
\]

\item \[
\frac{\partial^2\frac{f_\theta}{f}(z)}{\partial s_i\partial s_{i'}}(\Phi^0_t)=\left(e^2(z)-1\right)\phi({w^0_i}^Tx)\phi({w^0_{i'}}^Tx)
\]

\item  For $j\in\{t_{i-1}+1,\cdots,t_i\}$, let us write $\frac{\partial^2\frac{f_\theta}{f}(z)}{\partial s_i\partial w_j}:=\left(\frac{\partial^2 \frac{f_\theta}{f}(z)}{\partial s_i\partial w_{j0}},\cdots,\frac{\partial^2\frac{f_\theta}{f}(z)}{\partial s_i\partial w_{jd}}\right)^T$
\[
\frac{\partial^2\frac{f_\theta}{f}(z)}{\partial s_{i}\partial w_{j}}(\Phi^0_t)=\left(e^2(z)-1\right)\phi({w^0_{i}}^Tx)q_ja^0_i\phi^{'}({w^0_{i}}^Tx)x+e(z)q_j\phi^{'}({w^0_{i}}^Tx)x
\]
\item  For $j\in\{t_{i'-1}+1,\cdots,t_{i'}\}$, with $i\neq i'$, let us write $\frac{\partial^2\frac{f_\theta}{f}(z)}{\partial s_i\partial w_j}:=\left(\frac{\partial^2 \frac{f_\theta}{f}(z)}{\partial s_i\partial w_{j0}},\cdots,\frac{\partial^2\frac{f_\theta}{f}(z)}{\partial s_i\partial w_{jd}}\right)^T$
\[
\frac{\partial^2\frac{f_\theta}{f}(z)}{\partial s_{i}\partial w_{j}}(\Phi^0_t)=\left(e^2(z)-1\right)\phi({w^0_{i}}^Tx)q_ja^0_{i'}\phi^{'}({w^0_{i'}}^Tx)x
\]

\item For $j\in\{t_{k^0}+1,\cdots,k\}$ :
 \[
\frac{\partial^2\frac{f_\theta}{f}(z)}{\partial s_i\partial a_j}(\Phi^0_t)=\left(e^2(z)-1\right)\phi({w^0_i}^Tx)\phi({w_j}^Tx)
\]

\item For $j\in\{t_{i-1}+1,\cdots,t_i\}$  et $j'\in\{t_{i'-1}+1,\cdots,t_{i'}\}$, $j\neq j'$, let us write
\[\frac{\partial^2\frac{f_\theta}{f}(z)}{\partial w_j\partial w_{j'}}:=\left(
\begin{array}{ccc}
\frac{\partial^2 \frac{f_\theta}{f}(z)}{\partial w_{j0}\partial w_{j'0}}&\cdots&\frac{\partial^2\frac{f_\theta}{f}(z)}{\partial w_{j0}\partial w_{j'd}}\\
\vdots&\ddots&\vdots\\
\frac{\partial^2 \frac{f_\theta}{f}(z)}{\partial w_{jd}\partial w_{j'0}}&\cdots&\frac{\partial^2\frac{f_\theta}{f}(z)}{\partial w_{jd}\partial w_{j'd}}
\end{array}
\right)
\]
We get
\[
\frac{\partial^2\frac{f_\theta}{f}(z)}{\partial w_j\partial w_{j'}}(\Phi^0_t)=\left(e^2(z)-1\right)q_jq_{j'}a^0_i\phi^{'}({w^0_i}^Tx)a^0_{i'}\phi^{'}({w^0_{i'}}^Tx)xx^T
\]

\item For $j\in\{t_{i-1}+1,\cdots,t_i\}$, let us write 
\[\frac{\partial^2\frac{f_\theta}{f}(z)}{\partial w^2_j}:=\left(
\begin{array}{ccc}
\frac{\partial^2 \frac{f_\theta}{f}(z)}{\partial w^2_{j0}}&\cdots&\frac{\partial^2\frac{f_\theta}{f}(z)}{\partial w_{j0}\partial w_{jd}}\\
\vdots&\ddots&\vdots\\
\frac{\partial^2 \frac{f_\theta}{f}(z)}{\partial w_{jd}\partial w_{j0}}&\cdots&\frac{\partial^2\frac{f_\theta}{f}(z)}{\partial w^2_{jd}}
\end{array}
\right)
\]
We get
\[
\frac{\partial^2\frac{f_\theta}{f}(z)}{\partial w^2_j}(\Phi^0_t)=\left(e^2(z)-1\right)\left(q_ja^0_i\phi^{'}({w^0_i}^Tx)\right)^2xx^T+e(z)a^0_iq_j\phi^{''}({w^0_i}^Tx)xx^T
\]

\item For $j\in\{t_{i-1}+1,\cdots,t_i\}$ and $l\in\{t_{k^0}+1,\cdots,k\}$, let us write $\frac{\partial^2\frac{f_\theta}{f}(z)}{\partial w_j\partial a_l}:=\left(\frac{\partial^2 \frac{f_\theta}{f}(z)}{\partial w_{j0}\partial a_l},\cdots,\frac{\partial^2\frac{f_\theta}{f}(z)}{\partial w_{jd}\partial a_l}\right)^T$ :

\[
\frac{\partial^2\frac{f_\theta}{f}(z)}{\partial w_j\partial a_l}(\Phi^0_t)=\left(e^2(z)-1\right)q_ja^0_i\phi^{'}({w^0_i}^Tx)\phi({w_l}^Tx)x
\]

\item For $j,l\in\{t_{k^0}+1,\cdots,k\}$ :
 \[
\frac{\partial^2\frac{f_\theta}{f}(z)}{\partial a_j\partial a_l}(\Phi^0_t)=\left(e^2(z)-1\right)\phi({w_j}^Tx)\phi({w_l}^Tx)
\]

\end{itemize}
 Now, the terms:
\[
\left(1-\frac{1}{e^2(z)}\right)\left((\Phi_t-\Phi^0_t)^Tf^{'}_{(\Phi^0_t,\psi_t)}(z){f^{'}_{(\Phi^0_t,\psi_t)}}^T(z)(\Phi_t-\Phi^0_t)\right)
\]
and
\[
\sum_{i=1}^{k^0}\sum_{j=t_{i-1}+1}^{t_i}(q_jw_{j}-w^0_{i})^Txs_{i}\phi^{'}({w^0_i}^Tx)
\]
will be negligible compared to the first order term \( (\Phi_t-\Phi^0_t)^Tf^{'}_{(\Phi^0_t,\psi_t)}(z)\) when $\Phi_t\rightarrow \Phi^0_t$, even if the term of first order is of the same order or negligible compared to the terms of the second order: 
\[
e(z)\times\left(\sum_{i=1}^{k^0}\sum_{j=t_{i-1}+1}^{t_i}q_j(w_{j}-w^0_{i})^Txx^T(w_{j}-w^0_{i})a^0_i\phi^{''}({w^0_i}^Tx)\right).
\]
So, the development will be valid if for $\Phi_t\neq \Phi^0_t$, $z$ exists so that  
\[
\begin{array}{l}
e(z)\times\left(\beta-\beta^0+\sum_{i=1}^{k^0}s_i\phi({w^0_i}^Tx)\right.\\
+\sum_{i=1}^{k^0}\sum_{j=t_{i-1}+1}^{t_i}q_j\left(w_{j}-w^0_{i}\right)^Txa^0_i\phi^{'}({w^0_i}^Tx)\\
\left.\sum_{i=1}^{k^0}\sum_{j=t_{i-1}+1}^{t_i}q_j(w_{j}-w^0_{i})^Txx^T(w_{j}-w^0_{i})a^0_i\phi^{''}({w^0_i}^Tx)\right)\neq 0.\\
\end{array}
\] 
This inequality is guarantied by the assumption H-4.

These equations yield us the expression of $(\Phi_t-\Phi^0_t)^Tf^{''}_{(\Phi^0_t,\psi_t)}(z)(\Phi_t-\Phi^0_t)$
 $\blacksquare$  

\subsubsection*{References}

\small{
[1] Amari, S. \& Park, H. \& Ozeki, T. (2006) Singularities affect dynamics of learning in Neuromanifolds
{\it Neural computation }{\bf 18}: 1007-1065.

[2] Cottrell, M., Girard, B., Girard, Y., Mangeas, M. and Muller, C.  (1995).
Neural modeling for time series: a statistical stepwise method for weight elimination.
{\it IEEE Transaction on Neural Networks}, {\bf 6}, 1355--1364.

[2] Dacunha-Castelle, D. \& Gassiat, E. (1999) Testing the order of a model using locally conic parametrization: Population mixtures and stationary ARMA process. {\it The Annals of Statistics }{\bf 27}: 1178-1209.

[3] Fukumizu, K. (1996) A regularity condition of the information matrix of a multilayer perceptron network.
{\it Neural networks }{\bf 9} (5): 871-879.

[4] Fukumizu, K. (2003) Likelihood ratio of unidentifiable models and multilayer neural networks
{\it The Annals of Statistics }{\bf 31}: 833-851.

[5] Gassiat, E. \& Keribin, C. (2000) The likelihood ratio test for the number of components in a mixture with Markov regime.
{\it ESAIM Probability and statistics }{\bf 4}: 25-52.

[6] Gassiat, E. (2002) Likelihood ratio inequalities with applications to various mixtures.
{\it Annales de l'Institut Henri Poincar\'e }{\bf 38}: 897-906.

[7] Liu, X \& Shao, Y. (2003) asymptotics for likelihood ratio tests under loss of identifiability
{\it The Annals of Statistics }{\bf 31}: 807-832.

[8] Sussmann H. J. (1992) Uniqueness of the weights for minimal feed-forward nets with a given input-output map.
{\it Neural networks }{\bf 5} : 589-593.

[8] van der Vaart, A.W. (1998) {\it Asymptotic statistics} Cambridge: Cambridge university Press

[9] White, H. (1992) {\it Artificial Neural Networks: Approximation and Learning Theory.}
Oxford: Basil Blackwell.

}

\end{document}